\documentstyle[12pt]{article}
\def\hybrid{\topmargin 0pt      \oddsidemargin 0pt
        \headheight 0pt \headsep 0pt
        \textwidth 160true mm       
        \textheight 231true mm         
        \marginparwidth 0.0in
        \parskip 0pt plus 1pt   \jot = 1.5ex}
\catcode`\@=11
\def\marginnote#1{}
\newcount\hour
\newcount\minute
\newtoks\amorpm
\hour=\time\divide\hour by60
\minute=\time{\multiply\hour by60 \global\advance\minute by-\hour}
\edef\standardtime{{\ifnum\hour<12 \global\amorpm={am}%
        \else\global\amorpm={pm}\advance\hour by-12 \fi
        \ifnum\hour=0 \hour=12 \fi
        \number\hour:\ifnum\minute<10 0\fi\number\minute\the\amorpm}}
\edef\militarytime{\number\hour:\ifnum\minute<10 0\fi\number\minute}

\def\draftlabel#1{{\@bsphack\if@filesw {\let\thepage\relax
   \xdef\@gtempa{\write\@auxout{\string
      \newlabel{#1}{{\@currentlabel}{\thepage}}}}}\@gtempa
   \if@nobreak \ifvmode\nobreak\fi\fi\fi\@esphack}
        \gdef\@eqnlabel{#1}}
\def\@eqnlabel{}
\def\@vacuum{}
\def\draftmarginnote#1{\marginpar{\raggedright\scriptsize\tt#1}}
\def\draft{\oddsidemargin -.5truein
        \def\@oddfoot{\sl preliminary draft \hfil
        \rm\thepage\hfil\sl\today\quad\militarytime}
        \let\@evenfoot\@oddfoot \overfullrule 3pt
        \let\label=\draftlabel
        \let\marginnote=\draftmarginnote
   \def\@eqnnum{(\theequation)\rlap{\kern\marginparsep\tt\@eqnlabel}%
\global\let\@eqnlabel\@vacuum}  }
\relax
\hybrid
\def\cp{{\bf CP}^n} \def\gggg{{\bf g}}
\title{Some aspects of braided geometry:\\ differential calculus,
tangent space, gauge theory}
\author{P.~Akueson,\\
D.~Gurevich,\\
ISTV, Universit\'e de Valenciennes
59304 Valenciennes, France}

\begin{document}
\maketitle
\begin{abstract}
A new approach is suggested to quantum differential calculus on
certain quantum varieties.
 It consists in
replacing quantum de Rham complexes with differentials satisfying
Leibniz rule by those which are
in a sense close to Koszul complexes from \cite{G1}. We  also introduce the
tangent space  on a quantum
hyperboloid  equipped with  an action on the
quantum function space and define the notions of quantum
(pseudo)metric and
quantum connection (partially defined) on it.
All objects are considered from the viewpoint of
flatness of quantum deformations.
A problem of constructing a flatly deformed quantum gauge theory is
discussed as well.
\end{abstract}



\renewcommand{\theequation}{{\thesection}.{\arabic{equation}}}
\def\r#1{\mbox{(}\ref{#1}\mbox{)}}
\setcounter{equation}{0}
\def\ot{\otimes}
\def\uq{U_q({\bf g})}
\def\uqs{U_q(sl(n))}
\def\uqsl{U_q(sl(2))}
\def\aam{A^{(m)}}
\def\aan{A^{(n)}}
\def\vn{V^{\ot n}}
\def\iin{I^{(n)}}
\def\inn{I^{n}}
\def\de{\delta}
\def\De{\Delta}
\def\RR{\bf R}
\def\CC{\bf C}
\def\cp{{\bf CP}^n}
\def\Mat{\rm Mat}
\def\Sym{\rm Sym}
\def\ahqc{A^c_{\h,\,q}}
\def\Om{\Omega}
\def\aqc{{\cal A}^c_{ 0\, q}}
\def\ahqc{{\cal A}^c_{\h\, q}}
\def\ac{{\cal A}^c_{0\, 1}}
\def\h{{\hbar}}
\def\ah{{\cal A}_{\hbar}}
\def\vv{V^{\ot 2}}
\def\thq{T(H_q)}
\def\am{A_{\mu}}
\def\ami{A_{\mu}^i}
\def\gq{{\gggg}_q}
\def\ggq{{\gggg}_q^{\ot 2}}
\def\sss{\cal S}
\def\Ob{Ob\,}
\def\h{{\hbar}}
\def\a{\alpha}
\def\na{\nabla}
\def\ve{\varepsilon}
\def\ep{\epsilon}
\def\bea{\begin{eqnarray}}
\def\eea{\end{eqnarray}}
\def\beq{\begin{equation}}          \def\bn{\beq}
\def\eeq{\end{equation}}            \def\ed{\eeq}
\def\nn{\nonumber}
\def\oH{{\overline H}}
\def\aaaa{\cal A}
\def\aaa{\cal A}
\def\End{{\rm End\, }}
\def\uEnd{\underline{\rm End}\,}
\def\Ker{{\rm Ker\, }}
\def\root{{\rm root\, }}
\def\diag{{\rm diag\, }}
\def\Im{{\rm Im\, }}
\def\Id{{\rm Id }}
\def\ra{{\rm rank\, }}
\def\Id{{\rm id\, }}
\def\Vect{{\rm Vect\,}}
\def\Hom{{\rm Hom\,}}
\def\uHom{\underline{\rm Hom}\,}
\def\tr{{\rm tr}}
\def\det{{\rm det}}
\def\codet{{\rm codet}}
\def\dim{{\rm dim}\,}
\def\id{{\rm id}\,}
\def\span{{\rm span}\,}
\def\udim{\underline{\rm dim}\,}
\def\dem{{\rm det}^{-1}}
\def\Fun{{\rm Fun\,}}
\def\de{\delta}
\def\al{\alpha}
\renewcommand{\theequation}{{\thesection}.{\arabic{equation}}}
\setcounter{equation}{0}
\newtheorem{proposition}{Proposition}
\newtheorem{conjecture}{Conjecture}
\newtheorem{corollary}{Corollary}
\newtheorem{theorem}{Theorem}
\newtheorem{definition}{Definition}
\newtheorem{remark}{Remark}
\def\R{{\cal R}}
\def\n{{\bf n}}   \def\va{\varepsilon}
\def\so{s_{\omega_1+\omega_{n-1}}}
\section{Introduction}
In this paper we consider some problems which can be gathered
together under a general name
"braided (quantum, twisted or q-deformed) geometry". This type of
geometry has had a real gold rush since the
creation of the
quantum group (QG) theory. This phenomenon is motivated by a common
desire to generalize methods
of ordinary geometry for the needs of mathematical physics since, in
accordance with a widespread
opinion, the future of this discipline is connected with models which
are
 covariant w.r.t.
 special Hopf algebras rather than to ordinary transformation groups.

Nevertheless, it turned out that not all objects of the ordinary
geometry have their consistent
q-analogues. For example, all attempts initiated by L.Woronowicz
\cite{W1}, \cite{W2}
to develop a bicovariant differential calculus on quantum function
space $\Fun(SL_q(n))$
with two properties:
 flatness of deformation of the differential
algebra\footnote{Let us recall that a deformation
${\aaa}\to\ah$ where $\h$ is a formal parameter is called {\em flat}
if
$$1.\,\,\ah/\h\ah={\aaa}\qquad  2.\,\,\ah\quad {\em and}\quad
{\aaa}[[\h]]={\aaa}\ot k[[\h]]$$
are isomorphic  as $k[[\h]]$-modules (the tensor product is
complete in the $\h$-adic topology).
Here we  consider only the objects related to the famous
Drinfeld-Jimbo
QG $\uq$.
Nevertheless, some of them can be generalized to non-quasiclassical
Hecke symmetries,  i.e. solutions of the quantum Yang-Baxter equation
(QYBE) whose "symmetric" and "skew symmetric" algebras  possess
non-classical Poincar\'e series, (cf.
\cite{G1}).} and Leibniz rule for the corresponding differential have
failed.

Moreover, such a differential calculus does not exist. It was shown in
\cite{AAM} by considering the corresponding quasiclassical object,
namely the graded Poisson-Lie structure, which is an extension of the
Sklyanin-Drinfeld bracket to the differential algebra (cf. also the
last section of  \cite{Ar}).

The problem is that  a consistent q-deformation of the differential
calculus which is well defined
on Lie group  $GL(n)$ (more precisely, on the corresponding
matrix algebra $\Mat(n)$), cf.  \cite{T},
 is not compatible with the constrains resulting from the
equation $\det_q=1$ where $\det_q$ is
 the quantum
determinant. However, some q-deformed differential algebra equipped
with a differential without
Leibniz rule exists in the $SL(n)$ case, cf. \cite{Ar}, \cite{FP1}.
The authors of \cite{FP1} recall a claim of L.Faddeev that the Leibniz
rule is not reasonable in the quantum case.

One of the main purposes of this paper is to suggest a regular
way to construct
de Rham type complexes without any Leibniz rule.
The essence of such a complex is close to that of
Koszul complex (of the first kind) introduced in \cite{G1} (we recall
the construction of such a complex in section 2).
In order to introduce a differential we
 fix a base in a q-deformed differential algebra and define it
only in this base. This saves us from checking the fact that the
differential respects the relations
which define the algebra in question. Moreover, we realize in the
classical case a spectral analysis
 of de Rham complex, i.e. we study the behavior of the classical
 differential on irreducible
components of the initial complex and define a quantum differential
with similar properties but in the
q-deformed category. This approach is realizable when the spectral
structure of the complex
in question is simple enough (it can be also applied to a non-quasiclassical case).

We apply this approach to a quantum hyperboloid.
By construction, its cohomology is just the same as in the classical
case.
(In general, the following conjecture seems to be very plausible:
 once a quantum de Rham complex is constructed in a proper way
it has for a generic $q$ 
the same cohomology as its classical counterpart, cf. for
example \cite{FP2}, \cite{HS}
for an illustration  of this conjecture.)
The corresponding construction is described in section 3.

In this connection we  also discuss the following
problem: what is a proper definition of the tangent space on the
quantum hyperboloid
(in other words, what is the
phase space corresponding to the quantum hyperboloid considered as a
configuration space)?

We introduce such a tangent space $\thq$ (which is treated as
$A$-module where $A$
is the quantum function space
in question) and equip it with  an action
$$\thq\ot A\to A$$
converting elements of the tangent space into "braided vector fields".
Let us remark that
our construction of  braided vector fields is realized 
without (once more!)
any Leibniz rule (cf. \cite{A}).  

We also introduce (in section 4)  the notions of a (pseudo)metric
and a connection (partially defined) on the tangent
space  on  the quantum hyperboloid. In all
our constructions we impose
only  two properties on any q-deformed object in question:
$\uqsl$-covariance and flatness of the deformation.

In  section 5 we consider the problem of
constructing a quantum
gauge theory from this viewpoint.
In spite of numerous attempts to generalize the classical gauge theory from the above viewpoint
 up to now this has not been satisfactory.
We are rather sceptic about a possibility to introduce
 a consistent q-deformation of the classical gauge theory.
We disscus this in section 5.

Throughout the whole paper the basic field $k$ is $\RR$ or
$\CC$ and the parameter $q\in k$
is assumed to be generic.

\section{De Rham and Koszul complexes: comparative description}

First, let us consider  some complexes related to the QG $\uqs$.
The most popular complexes of such a type are de Rham complexes
connected with the
first fundamental modules of the QG $\uqs$ \cite{WZ} and those
defined on the q-deformed
matrix algebra $\Mat(n)$ \cite{T}. 
Whereas the former one is
one-sided $\uqs$-covariant, the latter one is bicovariant.
(Such  complexes exist for any Hecke symmetry, see footnote 1.)

A de Rham complex related to the first fundamental
$U_q(SO(n))$-module was constructed in \cite{CSW}.

However, all the above complexes are, in a sense, objects of  quantum
(braided or q-) linear algebra than  of quantum geometry.
Quantum geometry deals rather with
quantum varieties different from  vector spaces. A typical example of
such a variety is $SL_q(n)$ defined by the equation $\det_q=1$
 (by an abuse of the language we speak about a variety although in
 fact we deal
with the corresponding "quantum function space"). As we
said above,
the quantum differential calculus well defined on the vector space
$\Mat(n)$ cannot be restricted to the variety in question if
we want it to be a flat deformation
of its classical counterpart and its differential to obey the
Leibniz rule.
We refer the reader to the survey \cite{I} where this problem
is discussed.

Another interesting class of varieties connected with the QG $\uq$
are quantum homogeneous spaces which are
one-sided $\uq$-modules (the products in the corresponding algebras are assumed to 
be $\uq$-covariant  in the following
sense
$$Z(a\cdot b)=Z_{(1)}a\cdot Z_{(2)}b,\,\,Z\in \uq,\,\, 
Z_{(1)}\ot Z_{(2)}=\De (Z)).$$
 A quantum homogeneous space is usually introduced  via
a couple of QG in the spirit of a homogeneous space 
$G/H$. However, it is desirable to have its more
explicit description
by some system of equations.  

An attempt to find such a system for
certain q-deformed $SL(n)$-orbits
in $sl(n)^*$ featured in \cite{DGK} where a two parameter family of quantum algebras
was constructed. The problem of which algebra of this family could be
considered as a q-analogue
of a commutative algebra was not so evident.
 In the following we consider a particular case
of these q-deformed orbits, namely that related to $\uqsl$ and called
{\em quantum hyperboloid}.
 Being equipped with a proper involution it becomes
Podles' {\em quantum sphere} \cite{P1} (more precisely a particular
case of Podles' quantum sphere
which is simply the "q-commutative" case; note that in this low-dimensional
case there is no problem with
understanding "q-commutativity").

The first attempt to construct a q-deformed differential calculus on a quantum sphere
was undertaken in \cite{P2}. However, the corresponding differential algebra
is not a flat deformation of its classical counterpart. In  section 3 we will
present another approach to introducing a quantum de Rham complex with the flatness property.

Let us evoke  now  another type of complexes connected (in particular) with the QYBE,
 namely Koszul complexes (we are still working in the framework
 of quantum linear algebra). Let $V$ be a vector
space over $k$ and $I\subset\vv$ be a  subspace of $\vv$.
Let us set
$$I^{(0)}=k,\,\,I^{(1)}=V,\,\,\iin =I\ot V^{\ot (n-2)}\cap V\ot I\ot
V^{\ot (n-3)}\cap...
\cap V^{\ot (n-2)}\ot I,\,\,n\geq 2$$
and consider the quadratic algebra
$$A=T(V)/\{I\}, \quad {\rm where}\quad \{I\}\qquad {\rm is\,\,
the \,\, ideal\,\,
generated \,\,by}\quad I$$
and $T(V)$ stands for the free tensor algebra of the space $V$.

Let $\aan$ be its homogeneous component of degree $n$. Note that
$A^{(0)}=k,\,\,A^{(1)}=V$ and
 $\aan,\,\,n\geq 2$ can be treated as the quotient
$$\vn/\inn \qquad {\rm where}\qquad \inn= I\ot V^{\ot (n-2)}+ V\ot
I\ot V^{\ot (n-3)}+...
+V^{\ot (n-2)}\ot I.  $$
Then the corresponding  Koszul complex is defined by
\beq
d:A\ot \iin\to A\ot I^{(n-1)},\qquad d(a\ot x\ot y)=ax\ot y\quad
{\rm  where}\quad a\in A,\,\, x\ot y\in V\ot V^{\ot (n-1)}\label{com}\eeq
and $ax$ is the product in the algebra $A$.
In fact, this complex decomposes into a series of subcomplexes
$$A^{(m)}\ot \iin \to A^{(m+1)}\ot I^{(n-1)}.$$

\begin{definition} A quadratic algebra $A$ is called Koszul if the
cohomology of the complex \r{com}
vanishes  in all terms (except of course
the trivial term $A^{(0)}\ot I^{(0)}$ consisting of constants, i.e.
elements of $k$).
\end{definition}

Let us suppose now that we have two nontrivial complementary
subspaces $I_+\subset \vv$ and $I_-\subset \vv$, i.e. such that
$I_+\cap I_-=\emptyset$ and $I_+\oplus I_-=\vv$
 and associate to them two algebras
$$A_+=T(V)/\{I_-\}\qquad {\rm and}\qquad A_-=T(V)/\{I_+\}$$
(they are treated as "symmetric" and "skew symmetric" algebras
whereas the elements of the subspaces
$I_{\pm}\subset\vv$ are treated as "symmetric" and "skew symmetric"
tensors). Then we can define
two Koszul complexes
$$d:\iin_+\ot A_-\to  I^{(n-1)}_+\ot A_- \qquad {\rm and}  \qquad
\de :A_+\ot \iin_-\to A_+\ot  I^{(n-1)}_-$$
as it is described above.

If moreover, we can identify
$A^{(n)}_+$ with $\iin_+$ and $A^{(n)}_-$ with
$\iin_-$ (this means that the spaces $\iin_+$ and $I^n_-$ on the one
hand and $\iin_-$
and $I^n_+$ on the other hand are complementary for $n\geq 3$, cf.
\cite{DS})
 we can consider these two complexes as one (whose the terms are
$A_+^{(m)}\ot A_-^{(n)}$)
but equipped with two
 differentials
mapping in opposite directions.

This is just the case of the complexes constructed in \cite{G1} (where $V$ is a 
vector space equipped with a Hecke symmetry, cf. footnote 1).
 As  shown in \cite{G1}, the algebras $A_{\pm}$  are
Koszul. In particular, this implies  the classical relation
$$P_+(t)\, P_-(-t)=1$$ between the Poincar\'e series of the
"symmetric" and "skew symmetric" algebras.

Let us remark that the above identification $ A_{\pm}^{(n)}\approx
\iin_{\pm}$ was realized in \cite{G1} by means of projectors
$$P^n_{\pm}: V^{\ot n}\to \iin_{\pm}$$
whose kernels are just $\inn_{\mp}$. This implies that the spaces
$\iin_{\pm}$ and
$\inn_{\mp}$ are complementary. Moreover, the differentials $d$ and
$\de$ are realized in
\cite{G1} directly in terms of these projectors.

We say that an element $x\ot y \in A_+^{(m)}\ot A_-^{(n)}$ is given
in a canonical
(or base) form if $x\ot y\in I_+^{(m)}\ot I_-^{(n)}$, i.e. it is
realized as a sum of products of
"symmetrized" and "skew symmetrized" elements.  In virtue of \cite{G1}
 any element of $A_+^{(m)}\ot A_-^{(n)}$ can be represented in a canonical form.

Let us now compare these complexes with the de Rham complex
constructed in \cite{WZ}.
By applying the de Rham differential to the product
$x_{i_1}x_{i_2}...x_{i_{m}}$
one obtains, by virtue of the Leibniz rule, a sum whose arbitrary summand
is of the  form
$$x_{i_1}x_{i_2}...x_{i_{p-1}}dx_{i_p}x_{i_{p+1}}...x_{i_{m}},\quad
1\leq p\leq m.$$
(The sign $\ot$ is systematically omitted.) Here $\{x_i\}$ is a base of the space $V$
equipped with a Yang-Baxter operator of the Hecke type.

The second step of the procedure consists of moving  the factor
$dx_{i_p}$ to the
right  side (for concreteness). So, the problem arises of finding a
moving which would be compatible with
the differential and would lead to a flat deformation of the initial
differential algebra.
If,  moreover, one wants to restrict the differential
to a quantum variety it is necessary to coordinate such a movement with
constrains arising from the system of equations defining the variety
in question.

Nevertheless, such a problem does not
appear for the Koszul complex \r{com} since its
differential $d$ takes only one (namely, extreme) factor of the
space $\iin$ to the
algebra $A$. So, one should not transpose the elements from $V$ and
their differentials.

We can say that the Koszul complex from \cite{G1} and the de Rham one from
\cite{WZ} are formed by the same terms. The difference is that all elements of the
Koszul complex are
represented in  the canonical form. Moreover, it is easy to see that
the differentials of
these two complexes are
proportional to each other on each term (and the coefficients are not trivial).
This implies that their cohomologies are isomorphic
 (recall that $q$ is generic).

Since the cohomology of the Koszul complex is trivial (apart from the
$(0,0)$  term) we obtain
that it is also true for  de Rham complex from
\cite{WZ} (a quantum version of the Poincar\'e lemma).

Let us remark that this scheme can be extended to other couples of
subspaces $I_{\pm}$  associated to the QYBE (including non-quasiclassical cases)
but the crucial problem is to show that the associated spaces
$\iin_+$ and $I^n_-$
(resp., $\iin_-$ and $I^n_+$) are complementary (cf. \cite{DS}).

\section{De Rham type complex on quantum hyperboloid}

Let us pass now to a quantum hyperboloid. 
Consider the QG
$\uqsl$ generated by the generators $X,\,H,\,Y$ subject to the
well known relations
(cf. \cite{CP}). 
Let us fix a coproduct and the corresponding
antipode and consider the spin 1
$\uqsl$-module $V=V^q$. 

In order to define a quantum hyperboloid we should fix
 a base in $V$ and write down the
system of relations
on the generators compatible with action of the QG in question.
However, we want to represent this system in a symbolic way without
referring to its specific coordinate form.

We need only the fact that the fusion ring for $\uqsl$-modules is
exactly the same as in the classical
case (we consider only the finite-dimensional
$\uqsl$-modules which are deformations of the 
$sl(2)$-modules).
Thus, if $V_i$ is the spin $i$
 $\uqsl$-module then the classical formula
$$V_i\ot V_j=\oplus_{k=\vert i-j\vert}^{i+j} V_k$$
is still valid although the Clebsch-Gordan coefficients (which depend
on a  base)
are q-deformed.

In particular, we have
$$\vv=V_0\oplus V_1\oplus V_2.$$
We keep the notation $V$ for the initial space and $V_1$ for the
component in $\vv$ isomorphic
to $V$. Let us fix in the spaces $V,\,\,V_0,\,\, V_1$ and $V_2$ some
highest weight (h.w.) elements $v,\,\,v_0,\,\, v_1$ and $v_2$
respectively, and impose the relations
(which are the most general relations compatible with action of the
QG $\uqsl$)
\beq
v_0=c,\,\,v_1=\h v\;. \label{qh}
\eeq
Here $c\in k$ and $\h\in k$ are some constants.
One can now deduce the complete system of equations by applying to
the second relation the
decreasing operator $Y\in\uqsl$.

Let us denote $\ahqc$ the algebra defined by \r{qh} and derivative relations.

This algebra possesses the following property:
it is multiplicity free. More precisely,  any integer
spin module occurs once in its decomposition into a direct sum of
irreducible $\uqsl$-modules.
 Moreover, any element of $\ahqc$ can be represented in a unique way
as a sum of homogeneous elements
belonging to the components $V_i\subset V^{\ot i}$ (a proof of this fact can 
be deduced, for example from \cite{GV}). This
representation will be called {\em canonical or
base}. Note that element $v^{\ot i}$ is a h.w. one of the component
$V_i$.

We treat a particular case of the algebra in question, namely
$\aqc$, as a q-analogue of a commutative algebra 
 and call it {\em quantum
hyperboloid} if $c\not=0$ and
{\em quantum cone} if $c=0$.
 Since the  $sl(2)$-module $sl(2)^{\ot 2}$ is
multiplicity free  we can introduce
q-analogues $I_{\pm}$ of symmetric an skew symmetric subspaces of
$sl(2)^{\ot 2}$
by setting similarly to the classical case
$$I_+=V_0\oplus V_2 \quad {\em and } \quad I_-=V_1.$$

Let us emphasize that the corresponding algebras
$A_{\pm}=T(sl(2))/\{I_{\mp}\}$ are flat deformations of their
classical counterparts.

We will need also a q-deformed (braided) Lie bracket. It can be defined  as a non-trivial map
$$[\,\,,\,\,]_q:\vv\to V,$$
($V=sl(2)$ as linear spaces) being a $\uqsl$-morphism. By this request the bracket is defined 
in a unique way up to a factor. 

\begin{remark} Let us remark that for the  Lie algebras
$\gggg=sl(n),\, n\geq 3$ the $\gggg$-module
$\gggg^{\ot 2}$ is not multiplicity free any more: it possesses two
components
isomorphic to $\gggg$ itself, one belongs to the symmetric part of
$\gggg^{\ot 2}$
and the other one to the skew symmetric part.
 This is reason why it is no so evident what are q-analogues of the
 symmetric
and skew symmetric algebras of the space $\gggg$. However, there
exists a subspace $I_-\subset \ggq$
where $\gq=sl(n)$  as vector spaces but equipped with
a $\uqs$-module structure such that
the quadratic algebra
$T(\gq)/\{I_-\}$ is a flat deformation of the symmetric algebra
of $\gggg$ (cf. \cite{D}).
A more explicit description of $I_-$
 can be given by means of the so-called reflection equation (RE)
\beq
S\,L_1\,S\,L_1-L_1\,S\,L_1 S=0  \label{re} \eeq
where $S$ is a solution of the QYBE (here of the Hecke type),
 $L_1=L\ot \id $ and  $L$ is a matrix with matrix elements $(l_i^j)$.
The quadratic algebra defined by the system \r{re} is usually
called RE algebra.

Let us remark that the RE algebra is covariant w.r.t. 
$\uqs$ (cf. \cite{IP}) and it is a flat deformation  of its classical 
counterpart $\Sym(W)$ where $W=\span (l_i^j)$ (cf. \cite{L}).
It is not difficult to see that 
 the space $W$ is a sum of
two irreducible $\uq$-modules: one-dimensional one
with a generator $l=\tr_q\,L$ where $\tr_q$ is the q-trace
and $n^2-1$-dimensional one which can
be identified with $\gq$ above. By killing the component $l$ 
(i.e. by passing to the quotient of the RE algebra over the ideal
$\{l\}$) we get exactly the algebra mentioned above $T(\gq)/\{I_-\}$.
In other words, the space $I_-$ is defined by the relation \r{re}
but with one component less.

Moreover, by means of the RE algebra one can get an algebra looking like the enveloping algebra
 of q-deformed
 Lie algebra $sl(n)$ introduced in \cite{LS}.
Before killing the component $l$ let us realize a shift  $l_i^j\to l_i^j+\h\de_i^j$.
 Then instead of a graded quadratic algebra we get a filtered
 algebra, defined by quadratic-linear relations. Now by killing $l$ we get a
quadratic-linear algebra with $n^2-1$
generators. This is just
another realization of the enveloping algebra from \cite{LS} 
(if in the latter algebra we replace the Casimir element by a constant, cf. \cite{LS})
and  a flat two parameter deformation
of $\Sym(\gggg)$ whose  existence was
stated in \cite{D}. (However, to get a reasonable quasiclassical limit
we should replace the parameter $\h$ in this 
quadratic-linear algebra by $\h/(q-1)$.)

If $\gggg$ is a  simple Lie algebra different from $sl(n)$ its
tensor square is multiplicity free.
This allows one  to define a q-deformed Lie bracket requiring it
to be a non-trivial morphism
in the category of $\uq$-modules (this  defines the bracket uniquely
up to a factor)
and to introduce the enveloping algebra of the corresponding 
"braided Lie algebra" $\gq$.
Deformed analogues of the symmetric and skew symmetric algebras of the space $\gggg$ are also
well defined.
However, these algebras are not flat deformations of their classical
counterparts (cf. \cite{G2}).
\end{remark}

\begin{remark}
Let us emphasize that the classical counterpart $\ac$ of the algebra
$\aqc$ contains only
polynomials restricted to the hyperboloid (or the cone). This is
the reason why its properties and those of the function algebra on the sphere are similar : 
a passage from one algebra to the other one can be realized by a change of base. For
example, this passage does not change the cohomology of the 
de Rham complex (see below).
\end{remark}

Let us set $\Om^0=\aqc$.
Our next step is to define
the spaces of first- and second-order differential forms  over this
algebra.
First, consider the tensor products
$$\wedge^1=\aqc\ot V'\quad {\rm and} \quad \wedge^2=\aqc\ot V_1''.$$
Their second factors are treated as  pure
differentials (say
the element $x\ot y\in \wedge^1$
is treated as $x\,dy$ and in $x\ot y\in \wedge^2$ the factor $y\in
V_1''$ is a sum of  products
 of two pure first order differentials). The mark $'$  stands for a
 pure first-order
differential term and that $''$ stands for a pure second order
differential term. Thus, the space $V'$ (resp., $V_1''$) is isomorphic to the space 
$V_1$ itself; the isomorphism is defined by
$$d\,x_i\to x_i\,\,({\rm resp.,}\,\, d\,x_i\ot d\,x_j\to x_i\ot x_j).$$

Note that we treat the vector spaces $\wedge^i$ as left
$\aqc$-modules.
We do not endow their sum $\oplus \wedge^i$ with any algebraic
structure. So, we do not need any
transposition rule for the elements of $\aqc$ and their
differentials.

Let us introduce  now the  first- and second-order differential forms
on the quantum hyperboloid by
$$\Om^1=\wedge^1/\{(V\ot V')_0\},\,\,\Om^2=\wedge^2 /\{(V\ot
V_1'')_1+(v_0-c)\ot V_1''\}.$$
Here the terms in the denominators are  not ideals but only  left
$\aqc$-submodules of $\wedge^1$
and $\wedge^2$ respectively. The notation $(V\ot V')_i$ means that in
the
product $V\ot V'$ we
take the spin $i$ component (similarly for $(V\ot V_1'')_i$). And
$(v_0-c)\ot V_1'' $ stands for
the second-order differential forms containing $v_0-c$ as a factor.

To make this construction more explicit let us represent it in a
base form
(by restricting ourselves to the classical case since it does not
matter what case,
classical or quantum we deal with).
Let $u,v,w$ be the usual base in $\Fun(sl(2)^*)=\Sym(sl(2))$. Then the above
 denominators are  generated respectively by
$$2u\,dw+2w\,du+v\,dv \quad {\rm and} $$
$$2u e_2+v e_1,\,\,u e_3-w e_1,\,\,2w e_2+ v e_3,\,\,
(2uw+2wu+vv-c)e_i,\,\,i=1,2,3$$
$${\rm with} \quad
e_1=dudv-dvdu,\,\,e_2=dudw-dwdu,\,\,e_3=dvdw-dwdv.$$
(we suppose here that $v_0=2uw+2wu+vv$).

We have defined the spaces $\Om^i,\,\, i=1,\,2$ as some quotients.
Now, we want to
define the differentials in some bases of these spaces similarly 
to the Koszul complexes discussed above. To define such bases 
 we will realize a
spectral analysis of the spaces $\Om^i$, i.e. decompose these spaces into a direct sum of 
irreducible $sl(2)$-modules. First,  describe  the components in the products
$$V\ot  V',\,\,  V_i\ot V',\,\, V_i\subset\ac,\, i=2,\, 3,...$$
which are surviving in the quotient space $\Om^1$.
It is evident that in the product $V\ot V'$ only the components
$ (V\ot V')_1$ and $ (V\ot V')_2$ survive since by construction  the component
$ (V\ot V')_0$ is equal to 0 in the quotient.

By a similar reason in the product $V_2\ot V'$ the components
$(V_2\ot V')_2$ and $(V_2\ot V')_3$ survive and that $(V_2\ot V')_1$
is equal to 0 modulo
the terms of $k\ot V'=V'$. This can be explained as follows. The
elements
of $\mu^{12}(V\ot (V\ot V')_0)$ are trivial in $\Om^1$ by
construction.
Here $\mu$ stands for the product in $\aqc$, the indexes 12 mean as
usual that the operator $\mu$
is applied
to the first two factors. By reducing any element of the product
$V\ot V$
to the canonical form
we get a sum of an element from $V_2\subset \vv$ and another one from
$k$. This completes the proof.

Similarly,  in the product $V_i\ot V'$ the component
 $(V_i\ot V')_{i-1}$
is equal to 0 modulo the terms belonging to $V_j\ot V',\,\, j<i$.
Thus, we have shown the following.

\begin{proposition} The base in the $\aqc$-module $\Om^1$ is formed
by
$$1.\, V',\,\,2.\, (V\ot V')_{1, 2},\,\,3. \,(V_2\ot V')_{2, 3},\,\,
4.\, (V_3\ot V')_{3, 4}\quad
{\em etc}.$$
\end{proposition}

In a similar way one can perform a spectral analysis of the
$\aqc$-module $\Om^2$ and describe its
base.

\begin{proposition} The base in the $\aqc$-module $\Om^2$ is formed
by
$$1.\, V_1'',\,\,2.\, (V\ot V_1'')_{0,2},\,\,3.\, (V_2\ot
V_1'')_{3},\,\,4.\,(V_3\ot V_1'')_{4}
\quad {\rm etc}.$$
\end{proposition}

An evident difference between the modules $\Om^1$ and $\Om^2$
consists in the following.
The module $\Om^2$ is defined as a quotient of $\wedge^2$ over the sum of two
submodules.
Therefore,  two components in the products $V_i\ot
V_1'',\,\,i\geq 2 $
disappear and only one survives.
The component $V_1\ot V'$ is exceptional because the relation $v_0=c$
does not lead to any
constraint for  it.
Let us consider now the de Rham complex in the classical case
\beq
0\longrightarrow\Om^0\stackrel{d_0}\longrightarrow
\Om^1\stackrel{d_1}
\longrightarrow \Om^2\longrightarrow 0.\label{omm}
\eeq
Since the differential commutes with the $sl(2)$ action it takes any
irreducible $sl(2)$-module 
to either an isomorphic $sl(2)$-module or 0.
Using propositions 1 and  2 it is not difficult to describe the
irreducible $sl(2)$-modules of
$\Om^i,\,\,i=0,1,2$ belonging to
$\Ker\,\,d$ and those belonging to
$\Im\,\,d$.

\begin{proposition}
1. In $\Om^0$ the only  trivial  module, i.e. that consisting of the
elements of $k$ belongs to  $\Ker\,\,d_0$.
 $$2.\,\,\Ker\,d_1=V'\oplus (V\ot V')_2\oplus (V_2\ot V')_3\oplus
 (V_3\ot V')_4\oplus...$$
and therefore the modules
$$(V\ot V')_1,\,\, (V_2\ot V')_2,\,\,(V_3\ot V')_3,...$$
go to isomorphic modules in $\Om^2$.
\end{proposition}

\begin{corollary} The cohomology of the complex \r{omm} is the following one
$$\dim H^0=1,\,\, \dim H^1=0,\,\, \dim H^2=1,$$
$H^0$ is generated by 1 and $H^2$ is generated by $(V\ot V_1'')_0$.
\end{corollary}

Thus, it is just the cohomology of the sphere (see above, remark 2).

Let us extend now de Rham  complex \r{omm} to the quantum case.
The terms of
the quantum complex have just the same irredicible components as their
classical counterparts
 (but these components become $\uqsl$-modules).
 Now we should define differentials.
Define them on each term  by
requiring them to be 
$\uqsl$-morphisms and to be flat deformations of the classical
differentials (by this demand the differentials are defined  an each $\uqsl$-module
in a unique way up to a factor). Finally, we have by construction
 just the same cohomology as in the classical case.

 Comparing our construction with that from 
\cite{P2} we repeat that the latter one is not any flat deformation of 
its classical counterpart meanwhile our deformation is flat by construction. 
On the other hand, we have lost the structure of an algebra in the
$\aqc$-module $\Om=\oplus \Om^i$ and the Leibniz rule for the
differentials.

\section{Quantum tangent space and related structures}

In the present section we introduce the tangent space on
quantum hyperboloid and discuss some derived structures (metric,
connection).
We discuss also a way to realize
the q-deformed tangent space by means of "braided vector fields".
Hopefully, this approach is valid for other quantum varieties like
quantum orbits
considered in \cite{DGK}.
Similarly to the previous section we avoid using any specific
base form.

First, we consider a sphere $S^2$ given by $x^2+y^2+z^2=c$. Let $\Fun(S^2)$
be the space of  the polynomials restricted to
the sphere
and $\Vect (S^2)$ be the space of left vector fields. The latter space
is generated as a left $\Fun(S^2)$-module by three infinitesimal rotations
$$X=y\partial_z-z\partial_y,\,\,Y=z\partial_x-x\partial_z,\,\, Z=x\partial_y-y\partial_x.$$

It is easy to check that the vector fields $X,\,Y,\,Z$ satisfy the following relation
\beq
x\,X+y\,Y+z\,Z=0
\label{ts1}
\eeq
($x,\,y,\, z$ are treated here as operators via the product operator in the algebra $\Fun(S^2)$).
So, as a $\Fun(S^2)$-module $\Vect (S^2)$ can be realized as the quotient $M/N$ where 
$$M=\{a\,X+b\,Y+c\,Z,\,\, a,b,c\in \Fun(S^2)\},$$
$$N=\{f\,(x\,X+y\,Y+z\,Z),\,f\in \Fun(S^2)\}.$$

In what follows we call this $\Fun(S^2)$-module {\em tangent space} and denote it $T(S^2)$.
In fact it is just the space $\Vect (S^2)$ but we want to emphasize by this notation  
that we ignore the operator meaning of this space. 
As usual, the tangent space is introduced in local terms as a 
vector bundle. However, in the quantum case such a local description is not possible.

In a similar way there can be introduced the tangent space $T(H)$ on a hyperboloid $H$.
Namely, it  can  be realized as the quotient of a free $\ac$-module $M$ over its submodule 
$$N=\{f\,(2u\,W+2w\,U+v\,V)\}.$$
This is also motivated by the operator meaning of the generators: 
the generators $U,\,V,\,W$ are represented in the algebra $\Fun(H)=\ac$
by infinitesimal hyperbolic rotations.

Note that the symmetric algebra of the tangent space $T(S^2)$ 
 can be treated as the function algebra on the underlying 4-dimensional
algebraic variety embedded in the 6-dimensional space 
$$(\span (x,\,y,\,z,\,X,\,Y,\,Z))^* ; $$
this variety is defined by the equation of the sphere and that \r{ts1}
(if $k=\RR$ it is true for $c>0$).
A similar description is also valid for the symmetric algebra of $T(H)$.

Unfortunately, there does not exist any quantum analogue of this algebra
being its flat deformation (see below). Nevertheless, a reasonable q-deformation of tangent space 
equipped with an appropriated module structure exists.
The aim of this section is to describe this deformation, i.e., to introduce the tangent
space on the quantum hyperboloid as an $\aqc$-module and to realize its elements as operators 
looking like vector fields on the classical object. 

In order to do it we represent the defining relation of the tangent space $T(H)$ in 
a symbolic way:
\beq
(V\ot V')_0=0\label{tb}\eeq
 (hereafter the mark $'$ designs the space $\span(U,\,V,\,W)$). We treat the tangent space on the
hyperboloid as a left
$\ac$-module (as a right
$\ac$-module the tangent space can be given by $(V'\ot V)_0=0$).

It is evident that if we want
to define the tangent space on the quantum hyperboloid as a flat
deformation of its
classical counterpart we should use the same formula \r{tb} but in
the category of $\uqsl$-modules. Let us make a precise.
First, we introduce the left $\aqc$-module $\wedge^1$ as
in the previous section
but with another signification of the space $V'$. This means that the generators $du,\, dv,\,dw$
are replaced by $U,\,V,\,W$, while the $\uqsl$-module structure of $V'$ is unchanged. 
Second, we define the
{\em tangent space} on the quantum hyperboloid as its quotient like
$\Om^1$ above (fortunately, both the tangent and cotangent spaces as $\aqc$-modules are defined 
by the same equation \r{tb}).
Let us denote the quotient object by $\thq$ reserving the notation $H_q$ for the quantum
hyperboloid.

\begin{proposition} The $\aqc$-module $\thq$ is a flat deformation
of its classical counterpart.
\end{proposition}

Proof follows immediately from the explicit construction of the base
of this quotient
given in the previous section.

Let us assign now an operator meaning  to the elements of the  space $\thq$.

\begin{proposition} There exists a map 
\beq
\beta:\thq\ot\aqc\to\aqc.\label{anch}
\eeq
such that the diagram
$$
\begin{array}{ccc}
\aqc\ot\thq\ot\aqc&\longrightarrow  &\thq\ot\aqc\\
\downarrow&&\downarrow\\
\aqc\ot \aqc&\longrightarrow  &\aqc
\end{array}
$$
is associative. Here  the elements of $\aqc$ act on
$\aqc$ (the low arrow)
by the usual product. The vertical arrows are defined by means of $\beta$
and the top one makes use of the $\aqc$-module structure of $\thq$.
(Thus, the map $\beta$ realizes an action of the space $\thq$ on the algebra $\aqc$.)
\end{proposition}

This proposition allows us to realize the tangent space as an operator algebra where
the elements of the algebra $\aqc$ act  via the product operator.
We call the elements of the space $\thq$ (left) {\em braided vector fields} if the operators
$\beta(V')$   satisfy  the relations 
$$(\beta\ot\beta)(V'\ot V')_1-\sigma\,\beta\,[\,\,,\,\,]_q (V'\ot V')_1=0$$
where $[\,\,,\,\,]_q$ is the q-deformed Lie bracket introduced in section 3 and $\sigma\in k$ is a 
non-trivial factor. 
This means that $\beta$ realizes a representation (in the sense of  \cite{LS})
of the braided Lie algebra
defined by the bracket $\nu[\,\,,\,\,]_q$ with a proper factor $\nu$.

\begin{proposition} There exists a map $\beta$ from the previous proposition such that 
the elements of $\thq$ being represented  via $\beta$ becomes braided vector fields.
\end{proposition}

We refer the reader to \cite{A} for proofs of these statements
(the main idea of the construction has been suggested in \cite{DG2}). 
Here we only want  to say that the problem is to find good candidates for the
role of q-analogues of the infinitesimal hyperbolic rotations $U,V,W$.
They arise from the adjoint action of the q-Lie algebra $sl(2)_q$ onto itself
(note that the operators $X,H,Y$ coming from the QG $\uqsl$ do not
satisfy the relation \r{tb}).

Let us remark that similar statements  are valid for the tangent space
 treated as a right $\aqc$-module. 

Thus, we have an embedding
\beq
sl(2)_q \hookrightarrow \thq\label{qan}
\eeq
where  the tangent space is realized as braided vector fields space.
This embedding is a deformation of its classical counterpart
which is the simplest example of a so-called anchor (recall that an anchor consists of
an variety $M$, a Lie algebra $\gggg$ and an embedding of $\gggg$ into the vector field space
on $M$). This is reason why we call the embedding \r{qan} {\em quantum anchor} 
in spite  of the fact that the whole of the space $\thq$ is not equipped with any q-deformed
Lie bracket. We consider also the data $(\thq,\,\aqc)$ as a partial q-analogue of Lie-Rinehart
algebras \cite{R} ("partial" means here that the space $\thq$ is not equipped with any "q-Lie algebra"
structure properly coordinated with the product operator in the algebra $\aqc$).

After having represented the space $\thq$ by braided vector fields it is natural to  introduce
the space of braided differential operators as that generated by the braided vector fields and the
elements of $\aqc$ treated as 0-order operators (see above). In the classical case this space 
is spanned by the subspaces
$$\ac\ot {V'}^{\ot n}.$$
The fact that this algebra is closed w.r.t. the operator product
is assured by the Leibniz rule: by means of this rule it is possible to represent a product of two
elements of this form as a linear combination of such  elements.

Unfortunately, in the quantum case any form of the Leibniz rule does not exist
(this fact can be checked by direct calculations). 
Roughly speaking, this means that there does not exist any reasonable way to transpose the elements
of the algebra $\aqc$ and those of the space $V'$. This is also the reason why there does not exist
any "q-symmetric algebra" of the quantum tangent space $\thq$
being a flat deformation of its classical counterpart (see above).
Without going into detail we say only that the Yang-Baxter operator (arising from the 
universal R-matrix) being at first glance a good candidate for the role of such a  
transposition leads to a non-flat deformation  of the classical symmetric algebra.
(See also bellow, remark 3). 

Let us pass now to the problem of constructing a q-deformed metric on the
quantum tangent space.
To distinguish the quantum tangent spaces equipped with the left and
right $\aqc$-module structures
we will use for the first (second) one the notation $\thq_l$
($\thq_r$).

\begin{definition} We say that  an operator
$$<\,\,,\,\,> : \thq_l\ot_k\thq_r\to \aqc$$
is  quantum (pseudo-)metric if
it commutes with left and right multiplication by the elements from
$\aqc$ in the following sense
\beq
<fP, Q>=f<P,Q>, \, <P, Qf>=<P,Q>f,\,\, \forall
\,f\in\aqc,\,P\in\thq_l,
\,Q\in \thq_r\label{pair}
\eeq
(in particular, $P,\,Q\in V'$) and if it is compatible with the
action of $\uqsl$ .
The latter property means, as usual that
$$Z<\,\,,\,\,>= <\,\,,\,\,>\De(Z) ,\,\,\forall Z\in\uqsl$$
(this relation is treated as operator one in $\thq_l\ot_k\thq_r$).
A metric is called symmetric if
\beq
<\,\,,\,\,>(V'\ot V')_1=0.\label{sym}
\eeq
\end{definition}

\begin{proposition} There exists the unique (up to a factor)
symmetric quantum metric on the quantum
hyperboloid.
\end{proposition}

A proof of this fact is given in \cite{A}. We do not reproduce it
here. Let us indicate only
the crucial idea of the proof. First, it is necessary to describe
all pairings
$$<\,\,,\,\,>:V'\ot V'\to \aqc$$
compatible with the $\uqsl$ action. In order to do it we should
decompose
$V'\ot V'$ into a sum of the
irreducible $\uqsl$-modules. This gives rise to the following 
two parameter family of
$\uqsl$-covariant pairing
$$<\,\,,\,\,>(V'\ot V')_2= a V_2,\quad <\,\,,\,\,>(V'\ot V')_0=b$$
completed by relations \r{sym} (as usual, the relations are given in
a symbolic way). On the second step we should
impose
the condition
$$<\,\,,\,\,>^{23}(V\ot V')_0\ot V'=0$$
which results in a relation between the parameters $a$ and $b$. It
remains to verify that this
relation is compatible with the following one
$$<\,\,,\,\,>^{12}V'\ot (V'\ot V)_0=0$$
and then to extend the metric to the whole $\thq_l\ot_k\thq_r$ by
using the relations \r{pair}.

Let us emphasize that although we call the above pairing metric
("pseudo" means only that its classical
analogue is not positive definite) it is well defined on the product
of a left and a right
$\aqc$-modules. If we want now to define a similar pairing between
two left (or right)
$\aqc$-modules we should proceed in the following way.
Let us identify the left tangent space $\thq_l$ and the right one
$\thq_r$, i.e. define a map
$$\al: \thq_l\to\thq_r$$
being a $\uqsl$-isomorphism. 

Then on setting by definition
$$<X,Y>=<X,\al(Y)>,\,\, X,Y \in\thq_l$$
we get a pairing between two left $\uqsl$-modules. At first glance the map $\al$
can be defined by means of the YB operator
arising from the QG $\uqsl$. However, this operator which establishes a bijectivity between
free $\aqc$-modules $\aqc\ot V'$ and $V'\ot \aqc$ is not any
bijectivity on their factors $\thq_l$ and $\thq_r$ since it does not take 
the denominator corresponding to  $\thq_l$ to that corresponding to $\thq_r$. 
We suggest another way to define such a $\uqsl$-isomorphism $\al$.

Let us represent
the both objects as sums of $\uqsl$-modules in the spirit of
proposition 1. Then the $\uqsl$-morphisms
$$\al: (V\ot V')_1\to (V'\ot V)_1,\,\, (V\ot V')_2\to (V'\ot
V)_2,\,\,
(V_2\ot V')_2\to (V'\ot V_2)_2,...$$
are defined uniquely  up to a factor on each couple of components (for
the generating space $V'$ we put $\al=\Id$). 

However, we can reduce this freedom
by identifying the elements from
$$(V\ot V')_2\quad {\rm and} \quad (V'\ot V)_2,\qquad (V_i\ot
V')_{i+1}\quad {\rm and} \quad
(V'\ot V_i)_{i+1},\,\, i=2,\, 3,...$$
which coincide if we replace $V'$ by $V$. As for the components
$$(V\ot V')_1\quad {\rm and} \quad (V'\ot V)_1,\qquad (V_i\ot
V')_{i}\quad {\rm and} \quad
(V'\ot V_i)_{i},\,\, i=2,\, 3,...$$
their elements are identified
if  this operation leads to opposite images. It is not difficult to
see that in the classical case
this identification and that defined by the flip coincide (it is the
motivation of our method).

\begin{remark} Let us remark that for algebras looking like that
$\aqc$ but connected to an involutory YB operator an identification of their
left and right modules can be realized by means of this operator.
Non-involutivity of the YB operator arising from the QG $\uqsl$
which leads to the above defect prevents us also from a reasonable definition
of a tensor product $M_1\ot_{\aqc} M_2$ of two (say) left $\aqc$-modules. 
The problem is that there do not exist any reasonable way to transpose the factor
$f\in\aqc$ in the product
$$m_1\ot f\,m_2,\,\, m_1\in M_1,\,\, m_2\in M_2$$ 
on the left side so that the tensor product $\ot_{\aqc}$
is still associative and the module $M_1\ot_{\aqc} M_2$ is
a flat deformation of its classical counterpart assuming $M_1$ and $M_2$ to be flat deformations
of their classical counterparts. 
For an involutory YB operators this problem does not appear.
\end{remark}

Let us discuss now the problem of defining a (torsion free) $\uqsl$-covariant
connection on the tangent space $\thq$. Such a partially defined
connection was introduced in \cite{A}.
"Partially defined" means here that the operators of covariant derivatives are
 defined only on a subspace
of $\thq$, namely on $V'$. More precisely, there exists a $\uqsl$-morphism
$$\na:\,\thq_l\ot V'\to\thq_l,$$
$$X\ot Y\mapsto\na_XY$$
such that
$$\na_{fX} Y=f\na_X Y,\quad X\in \thq,\,\,Y\in V',\,\,f\in\aqc$$
and
\beq
a^{ij}\na_{X_i}X_j=[X_i,X_j]_q,\,\,X_i,\,X_j\in V'\label{br}
\eeq
where $a^{ij} X_i\ot X_j\in V_1$ and $[\,\,,\,\,]_q$ is the above
mentioned q-deformed Lie bracket.

We would be able to extend this (partially defined) connection to the whole
$\thq\ot\thq$ if we could extend the bracket $[\,\,,\,\,]_q$ to
the whole $\thq$ and to understand
what is the enveloping algebra of this extended q-deformed Lie
algebra (we need this in order to write suitable expressions
in the l.h.s. of \r{br}). Unfortunately, we
do not know any way to do it.

Let us remark that all naive extensions of the bracket $[\,\,,\,\,]_q$
is not compatible  with the equation \r{tb}.

There exists a number of papers introducing 
the notions of quantum metric and connection in another way (cf. \cite{HM} and the references
therein). Our approach
is motivated by our desire to control the flatness of deformation of classical objects
(see also the next section).

Completing this section we want to  mention a very important property
of the $\aqc$-module $\thq$  (if $c\not=0$): it is  projective  in the category of $\uqsl$-modules.
This means that it is a direct component in a free $\aqc$-module
and there exists a projector 
of the latter one onto the module in question being a $\uqsl$-morphism (cf. \cite{A}).

Some other projective modules over quantum sphere have been considered in \cite{HM}. We plan
to devote a subsequent paper to quantum projective modules in a more general context.

\section{On quantum gauge theory}

There exist two approaches to q-deformed gauge theory. One of them
deals with the
usual manifolds (varieties) and deforms  only a structure of fibers.
The second approach  deals with quantum varieties looking like
the quantum hyperboloid above.

First, let us evoke the paper \cite{S} as the most advanced
contribution to the first kind approach. 
The gauge potential $\am$ introduced in this paper is a  vector
field
$$\am(x)=\ami(x)\,X_i$$
with values in the quantum Lie algebra $\gq$. Here $\gq$ is the
q-deformed Lie algebra $\gggg=sl(n)$
(or  $su(n)$) as defined in \cite{LS} (note that the case
$n=2$ was previously considered in \cite{DG1}). Thus, the factors $X_i$ are
elements of this quantum Lie algebra
and those $\ami(x)$ are usual functions depending on a "space-time
point" $x$ (or more generally, on
a point of a usual variety).

In virtue of \cite{LS} the quantum Lie algebra $\gq$ is realized as a
subspace in $\uq$ so that it is stable w.r.t. the adjoint action of the QG $\uq$ on
itself and
$$\De(X_i)=X_i\ot C+u_i^j\ot X_j$$
where $\De$ is the coproduct in $\uq$, $C$ is a central (Casimir) element of
$\uq$ and
$u_i^j$ are some  elements of $\uq$.

The crucial point of any gauge theory is a transformation law of $\am$
under an action of a
gauge group or a gauge Hopf algebra. In \cite{S} it is supposed to be
\beq
\am\mapsto A_{\mu}'(x)=h(x)_{(1)}\am
s(h(x)_{(2)})-\al^{-1}s(C)^{-1}\partial_{\mu}(h(x)_{(1)})
s(h(x)_{(2)}) \label{gt}
\eeq
where $\al$ is a coupling constant, $s$ is the antipode, $h(x)$ is a
function of $x$ with values in
the QG  $\uq$ and
$$h(x)_{(1)}\ot h(x)_{(2)}=\De h(x).$$
 However,  a problem arises
to "distribute the $x$-dependence of the coproduct $\De(h(x))$
between two factors" in the second
term  of \r{gt} so that it becomes an element of  the space $\gq$
for a fixed $x$ (let us emphasize that it is not a trivial task).
 Such a distribution has been found in  \cite{S}. 

Nevertheless, it was indicated in \cite{S} that if we consider a
(say) bosonic field $\psi(x)$ then
its defining relations cannot be introduced in a way compatible with
quantum gauge transformations.
We will try to explain it as follows. Let us consider the "quantum
covariant derivative"
of the field $\psi$
$$D_{\mu}\psi=\partial_{\mu}\psi+\beta\rho(\am)\psi$$
where $\rho$ is the representation of the QG $\uq$
corresponding to $\psi$ and $\beta=\al\rho(s(C))$ is a constant.
 In the r.h.s. of this formula the operator $\partial_{\mu}$ commutes
 with the
$\uq$ action but the operator $\rho(\am)$ does not.  This implies
that such a covariant derivative
cannot preserve the relations valid for $\psi$. (We would have a
similar effect in a supertheory
if we allowed the summand $\am$ to be an odd operator.)

Let us discuss now the approach of the second kind, i.e. we
suppose that the base variety is quantum as well.
An axiomatic way to suggest such an approach was considered in
numerous papers. We do not give here
an exhausting list of these papers and only refer the reader to the
papers \cite{S} and \cite{BM} where such a list is given.

We will point out the common features of all of them. First, a
quantum variety in question is given
in a way which does not allow us to control the flatness of
deformation (as a rule this problem is not
even evoked). Another crucial defect of this approach is that a
connection is introduced
habitually via a Leibniz rule similarly to the classical case
but as we have seen 
this implied the non-flatness of the deformation.  Another reason of
the non-flatness of the deformation
is that in the formulae analogueical to \r{gt}
the second summand  does not belong usually to the fiber.

This explains our scepticism about a possibility of constructing a
quantum gauge theory
related to the QG $\uq$ which would be a flat deformation of its
classical counterpart.

Anyway, it would be desirable to precede any attempt to construct
such a theory by a quasiclassical study in the spirit of \cite{Ar}
confirming or refuting  a possibility to do it.

{\bf Acknowledgment} The authors are very indebted to P.Pyatov and P.Saponov
for discussions and valuable suggestions.
 One of the authors (D.G.) was supported by the grant
CNRS PICS-608.

\end{document}